\documentclass[11pt]{amsart}
\usepackage{palatino}
\usepackage[all]{xy}
\usepackage{amssymb}
\usepackage{latexsym}
\usepackage{amscd}
\usepackage{color}
\input amssym.def
\input amssym

\newtheorem{theorem}{Theorem}[section]
\newtheorem{lemma}[theorem]{Lemma}

\newtheorem{proposition}[theorem]{Proposition}
\newtheorem{corollary}[theorem]{Corollary}

\newtheorem{remark}{Remark}

\def \<{\langle}
\def \>{\rangle}

\def \a{\alpha }

\def \l{\lambda }
\def \ga{\gamma }

\def \b{\beta }

\newcommand{\bea}{\begin{eqnarray}}
\newcommand{\eea}{\end{eqnarray}}
\newcommand{\be}{\begin {equation}}
\newcommand{\ee}{\end{equation}}

\newcommand{\h}{\frak h}
\newcommand{\wt}{{\rm {wt} }   }

\newcommand{\hh}{\hat {\frak h} }

\newcommand{\Z}{\Bbb Z}
\newcommand{\lar}{\longrightarrow}
\newcommand{\W}{\mathcal W}
\newcommand{\mip}{\overline{M(1)_p}}
\newcommand{\Zp}{{\Bbb Z}_{>0} }

\newcommand{\N}{{\Bbb Z}_{\ge 0} }
\newcommand{\C}{\Bbb C}

\newcommand{\la}{\langle}
\newcommand{\ra}{\rangle}

\newcommand{\nb}{{\mbox {\tiny{ ${\bullet \atop \bullet}$}}}}
\newcommand{\nn}{\nonumber \\}
\usepackage{amssymb}

\begin{document}
\title[Logarithmic operators and $\mathcal{W}$-algebras]{Logarithmic intertwining operators
and $\mathcal{W}(2,2p-1)$-algebras}
\author{Dra\v{z}en Adamovi\'c and Antun Milas}

\begin{abstract}
\noindent For every $p \geq 2$, we obtained an explicit
construction of a family of $\mathcal{W}(2,2p-1)$-modules, which
decompose as direct sum of simple Virasoro algebra modules.
Furthermore, we classified all irreducible self-dual
$\mathcal{W}(2,2p-1)$-modules, we described their internal
structure, and computed their graded dimensions. In addition, we
constructed certain hidden logarithmic intertwining operators
among two ordinary and one logarithmic
$\mathcal{W}(2,2p-1)$-modules. This work, in particular, gives a
mathematically precise formulation and interpretation of what
physicists have been referring to as "logarithmic conformal field
theory" of central charge $c_{p,1}=1-\frac{6(p-1)^2}{p}, p \geq
2$. Our explicit construction can be easily applied for
computations of correlation functions. Techniques from this paper
can be used to study the triplet vertex operator algebra
$\mathcal{W}(2,(2p-1)^3)$ and other logarithmic models.
\end{abstract}

\address{Department of Mathematics, University of Zagreb, Croatia}
\email{adamovic@math.hr}

\address{Department of Mathematics and Statistics,
University at Albany (SUNY), Albany, NY 12222}
\email{amilas@math.albany.edu}

\maketitle

\renewcommand{\theequation}{\thesection.\arabic{equation}}
\renewcommand{\thetheorem}{\thesection.\arabic{theorem}}
\setcounter{equation}{0} \setcounter{theorem}{0}
\setcounter{equation}{0} \setcounter{theorem}{0}
\setcounter{section}{-1}

\section{Introduction}

The Virasoro algebra is one of the most fundamental structures in
two-dimensional conformal field theory. The most important family
of Virasoro algebra modules are certainly the minimal models,
because these models give rise to rational conformal field
theories. Interestingly, many non-rational models have recently
appeared in studies of $\mathcal{W}$-algebras, which are certain
extensions of Virasoro vertex algebras. Since there are several
different types of $\mathcal{W}$-algebras (see for instance
\cite{FKRW} for $\mathcal{W}$-algebras of positive integer central
charge), in this paper we limit ourselves to
$\mathcal{W}$-algebras closely related to representations of
Virasoro algebra with central charge
$$c_{p,1}=1-\frac{6(p-1)^2}{p}, \ p \in \mathbb{N}_{\geq 2}.$$
These central charges, belonging to the boundary of Kac's table,
are relevant in logarithmic conformal field theory \cite{F2},
\cite{F3}, \cite{Ga}, \cite{Gu}. If we denote by $L(c_{p,1},0)$
the simple lowest weight Virasoro algebra module of central charge
$c_{p,1}$ and conformal weight zero, then we have the following
embedding of $\mathcal{W}$-algebras
\begin{center}
$L(c_{p,1},0)$ $\hookrightarrow$  $\mathcal{W}(2,2p-1)$
$\hookrightarrow$ $\mathcal{W}(2,(2p-1)^3)$, \ $p \geq 2,$
\end{center}
where $\mathcal{W}(2,2p-1)$ is also known as the {\em singlet}
$\mathcal{W}$-algebra, and $\mathcal{W}(2,(2p-1)^3)$ is the {\em
triplet} $\mathcal{W}$-algebra. The theory of
$\mathcal{W}$-algebras could be understood much better from vertex
algebra point of view. In this setup, the singlet vertex algebra
$\mathcal{W}(2,2p-1)$ is generated by the Virasoro element
$\omega$ and another element $H$ of conformal weight $2p-1$ (cf.
\cite{A-2003}, \cite{EFH}, \cite{H},  \cite{KaW}). The singlet
vertex algebra admits infinitely many nonisomorphic irreducible
modules, so it fails to be rational. On the other hand, the
triplet algebra $\mathcal{W}(2,(2p-1)^3)$ (cf. \cite{Ka1},
\cite{Ka2}) is of the right size and its rationality was discussed
in \cite{GK1}, \cite{GK2} (some further studies of the triplet
algebra were pursued in  \cite{FST}, \cite{GR} etc.). The triplet
algebra has only finitely many equivalence classes of irreducible
modules, but in addition it admits certain indecomposable
logarithmic modules (i.e., modules that admit nontrivial Jordan
blocks with respect to the action of the degree zero Virasoro
generator). These indecomposable modules are needed to obtain the
{\em fusion closure} and the modular invariance of characters.
These logarithmic modules are also responsible for logarithmic
behavior of matrix coefficients. All these properties make the
triplet algebra CFT a rather odd looking "rational CFT"-so in
order to distinguish it from ordinary rational CFTs-the triplet
model and related models are usually dubbed as rational
logarithmic CFTs.

Many important aspects of logarithmic CFT can be studied by using
algebraic techniques. For instance, the appearance of
non-diagonalizable representations can be easily explained with
the use of Zhu's associative algebra \cite{Ga}, \cite{M2},
\cite{M1}. Similarly, logarithmic behavior of correlation
functions can be explained via logarithmic intertwining operators
 \cite{HLZ1}, \cite{HLZ2}, \cite{M2}, \cite{M1} (see also
\cite{FFHST} for a related approach). The most interesting
examples of logarithmic intertwining operators are those of type
\be { {\rm logarithmic \ module} \choose {\rm ordinary \ module} \
\ \ {\rm ordinary \ module}}. \nonumber \ee In \cite{M2} the
second author provided a general construction of intertwining
operators that arise from certain deformations of bosonic vertex
operators. These operators are not present in the original
non-logarithmic theory, which is the main reason why we called
them {\em hidden}.

Several clues from the physics literature indicate that the
triplet algebra and other related models should involve hidden
logarithmic intertwinners. In this paper we make the first step in
the direction of understanding these operators. Here we consider
the intermediate singlet vertex algebra $\mathcal{W}(2,2p-1)$,
which will be denoted by $\mip$ throughout the paper (cf.
\cite{A-2003}). We prove several general results about
$\mathcal{W}(2,2p-1)$-algebras. Firstly, we show that $\mip$ is a
simple vertex algebra (cf. Theorem \ref{simple}). Secondly, we
construct a distinguished irreducible $\mip$-module denoted by
$M(1,\beta)$, which decomposes further as a direct sum of simple
Virasoro algebra modules (cf. Theorem \ref{irr-module}). Then we
classify all irreducible self-dual $\mip$-modules and we compute
their graded dimensions (cf. Theorem \ref{irreducible-ga-i}).
Furthermore, we construct a self-extension of $M(1,\beta)$, which
give rise to an indecomposable $\mip$-module (see Theorem
\ref{self-dual-th}). Finally, by using this self-extension and
several results from \cite{M1} we obtain an explicit construction
of a family of hidden intertwining operators of $\mip$-modules
(cf. Corollary \ref{hidden}). Thus, we construct an algebraic
counterpart of logarithmic conformal field theory of level
$c_{p,1}$. In a sequel we plan to study the triplet model and
related models.

\section{Feigin-Fuchs modules}
We shall introduce some notation first. We denote by $\goth{h}$ a
one-dimensional abelian Lie algebra spanned by $h$ with a bilinear
form $\langle \cdot, \cdot \rangle$, such that $\langle h, h
\rangle =1$, and by
$$\widehat{\goth{h}}=\goth{h} \otimes \mathbb{C}[t,t^{-1}] +
\mathbb{C}c$$ the affinization of $\goth{h}$ with bracket
relations
$$[a(m),b(n)]=m \langle a, b \rangle \delta_{m+n,0}c, \ \ a,b \in \goth{h},$$
$$[c,a(m)]=0.$$
Set $ \hat{{\h}}^{+}=t{\C}[t]\otimes
{\h};\;\;\hat{{\h}}^{-}=t^{-1}{\C}[t^{-1}]\otimes {\h}. $ Then
$\hat{{\h}}^{+}$ and $\hat{{\h}}^{-}$ are abelian subalgebras of
$\hat{{\h}}$. Let $\mathcal{U}(\hat{{\h}}^{-})=S(\hat{{\h}}^{-})$
be the universal enveloping algebra of $\hat{{\h}}^{-}$. Let ${\l}
\in {\h}$. Consider the induced $\hat{{\h}}$-module
\begin{eqnarray*}
M(1,{\l})=\mathcal{U}(\hat{{\h}})\otimes
_{\mathcal{U}({\C}[t]\otimes {\h}\oplus {\C}c)}{\C}_{\l}\simeq
S(\hat{{\h}}^{-})\;\;\mbox{(linearly)},\end{eqnarray*} where
$t{\C}[t]\otimes {\h}$ acts trivially on ${\C}_{\l}={\C}$,
${h}$ acts as $\la h, {\l} \ra$, and
$c$ acts as multiplication by 1. For simplicity, we shall write
$M(1)$ for $M(1,0)$.

It is well-known that $M(1)$ has a vertex operator algebra
structure and that each $M(1,\lambda)$ is a $M(1,0)$-module
\cite{FrB}, \cite{LL}. More precisely, there are infinitely many
different (non-isomorphic) vertex operator algebra structures on
$M(1)$, denoted by $M(1)_a, a \in \mathbb{C}$, where the conformal
vector is chosen to be \be \label{vir-vector}
\omega_a=\frac{h(-1)^2{\bf 1}}{2}+a h(-2){\bf 1}. \ee Similarly,
each $M(1,\lambda)$, denoted now $M(1,\lambda)_a$, becomes an
irreducible $M(1)_a$-module. Here the subscript $a$ indicates that
Virasoro algebra acts differently; as vector spaces of course
$M(1,\lambda)=M(1,\lambda)_a$.

It is a standard fact (which can be easily shown) that the vertex
operator algebra $M(1)_a$ has central charge
$$c=1-12a^2.$$
Also, $M(1,\lambda)_a$ is a Virasoro algebra module of lowest
conformal weight
$$h_{\lambda}=\frac{1}{2}\lambda^2-\lambda a.$$
The Virasoro algebra modules $M(1,\lambda)_a$ are usually called
Feigin-Fuchs module \cite{FF}.
Let us ignore for a moment the conformal structure and view
$M(1,\lambda)$ only as a $\widehat{\goth{h}}$-module. If we denote
by $M(1,\lambda)^{\circ}$ the contragradient
$\hat{\goth{h}}$-module of $M(1,\lambda)$ defined by using the
anti-involution $\omega(h(n))=-h(-n)$, then we have
$$M(1,\lambda)^{\circ} \cong M(1,-\lambda).$$ But if we use
anti-involution $\omega(L(n))=L(-n)$ and denote by
$M(1,\lambda)^*_a$ the contragradient Virasoro module (or
$M(1)_a$-module) the result is different as illustrated by the
following lemma (cf. \cite{FF}). We will include the proof here
for later purposes (see Section 6).
\begin{lemma} \label{dual-iso}
We have the following isomorphism
$$M(1,\lambda)^*_a \cong M(1,2a-\lambda)_a,$$
of Virasoro algebra modules (or $M(1)_a$-modules). In particular,
if $a=\lambda$, then $M(1,\lambda)_{a}$ is self-dual.
\end{lemma}
\noindent {\em Proof.} From the formula for $L(-n)$ in terms of
Heisenberg algebra generators $h(n)$, we have \bea && \langle L(n)
\cdot w',w \rangle=\langle w', L(-n) w \rangle \nn && =\langle
w',\left(\sum_{n \in \mathbb{Z}} \nb h(-m)h(-n+m) \nb - (-n+1)a
h(-n) \right) w \rangle \nn && =\langle \left(\sum_{n \in
\mathbb{Z}} \nb h(m)h(n-m) \nb - (n-1)a h(n) \right)w',w \rangle
\nn && = \langle \left(\sum_{n \in \mathbb{Z}} \nb
\bar{h}(m)\bar{h}(n-m) \nb -  (n+1) a \bar{h}(n) \right)w',w
\rangle \nn && = \langle \bar{L}(n) w',w \rangle, \eea where $\nb
\ \nb$ denotes the normal ordering,
$$\bar{h}(n)=h(n)-2a \delta_{n,0},$$ for all $n \in \mathbb{Z}$
and $\bar{L}(-n)$ are Virasoro generators in terms of $\bar{h}(n)$
generator. The map $h(n) \mapsto \bar{h}(n)$ induces an
automorphism of $\widehat{\goth{h}}$, and with this new Heisenberg
algebra generators the module is isomorphic to the dual of
$M(1,\lambda-2a)$, which is isomorphic to $M(1,2a-\lambda)_a$.
\qed

The embedding structure of Feigin-Fuchs modules is well-known
\cite{FF}. In the self-dual case it is particularly simple.
\begin{proposition} \label{self-dual}
As before, we let
$$\lambda_p=\frac{p-1}{\sqrt{2p}}, \ \  p \in \mathbb{N}_{\geq 2}.$$
The Feigin-Fuchs module $M(1,a)_{a}$ is self-dual and completely
reducible if and only if $a=\lambda_p$. Moreover, we have the
following decomposition
$$M(1,\lambda_p)_{\lambda_p}= \bigoplus_{n = 0}^\infty L(c_{p,1},h^p_n)$$
where $L(c_{p,1},h_n^p)$ denotes irreducible lowest weight
Virasoro module of central charge $c_{p,1}$ and lowest conformal
weight
$$h_n^p=\frac{(2pn)^2-(p-1)^2}{4p}.$$
\end{proposition}
In the previous proposition $a=\lambda_p$, so the central charge
is \be \label{cc} c_{p,1}=1-\frac{6(p-1)^2}{p} \ee and the lowest
conformal weight of $M(1,\lambda_p)_{\lambda_p}$ is
$$h_0^p=-\frac{(p-1)^2}{4p}.$$
Modules of central charge $c_{p,1}$ are also known as {\em
logarithmic} minimal models in the physics literature (this should
not to be confused with logarithmic modules that will appear later
in the text).

\section{Virasoro Verma modules of central charge $c_{p,1}$}

In the previous section we considered some special Feigin-Fuchs
modules. Here we discuss closely related Verma modules. Their
embedding structure is similar to those of Feigin-Fuchs modules
(for a fixed $c$ and $h$ one has to "invert" one-half of
embeddings in the Verma module to get the embeddings in the
Feigin-Fuchs module with the same $c$ and $h$). As usual, we shall
denote by $V(c,h)$ the Verma module of lowest conformal weight $h$
and central charge $c$, i.e.,
$$V(c,h)=\mathcal{U}(Vir) \otimes_{\mathcal{U}(Vir)_{\geq 0}} \mathbb{C}v_{c,h},$$
where $L(n)$, $n \geq 1$ acts trivially on the lowest weight
vector $v_{c,h}$ and \bea L(0) \cdot v_{c,h}&=& h v_{c,h} \nn C
\cdot v_{c,h}&=&c v_{c,h}. \nonumber \eea The Verma module is a
cyclic $\mathbb{N}$-gradable Virasoro algebra module, where the
grading is inherited from the action of $L(0)$. There exists a
unique maximal submodule of $V(c,h)$, denoted by $V^1(c,h)$, not
necessarily cyclic, such that $L(c,h)=V(c,h)/V^1(c,h)$ is
irreducible. In fact, every irreducible lowest weight module of
central charge $c$ and lowest conformal weight $h$ is isomorphic
to $L(c,h)$.

Generically, Verma modules are irreducible and isomorphic to
appropriate Feigin-Fuchs modules described in the previous
section. For instance, for each $a$, $-a^2/2 \notin \mathbb{Q}$,
we have an isomorphism $M(1,a)_a \cong L(1-12a^2,-\frac{a^2}{2})$.
This class of representations will not be treated in our paper.

%
%
In addition to complete description of Feigin-Fuchs  modules,
Feigin and Fuchs classified all embeddings among Verma modules.
Here is their result in the case of $c_{p,1}$ (cf. \cite{FF}):
\begin{proposition} \label{VermaVir}
The Verma module $V(c_{p,1},h)$ is reducible if and only if
$$h=h_{m,n}:=\frac{(m-np)^2-(p-1)^2}{4p}, \ \ n=1, \ \ m>0 \footnote{Equivalently ,
we may assume that  $n >0, \ \ 0 < m \leq p$}.$$ Moreover, for
$m=k p$, $k \in \mathbb{N}$ and $n=1$ we have the following chain
of embeddings
$$V(c_{p,1},h_{m,1}) \longleftarrow V(c_{p,1},h_{m+2p,1}) \longleftarrow
V(c_{p,1},h_{m+4p,1}) \longleftarrow V(c_{p,1},h_{m+6p,1})
\longleftarrow \cdots,$$ while for $1 \leq m \leq p-1$ we have
$$V(c_{p,1},h_{m,1}) \longleftarrow V(c_{p,1},h_{m+2p,1})
\longleftarrow V(c_{p,1},h_{-m+4p,1}) \longleftarrow
V(c_{p,1},h_{m+4p,1}) \longleftarrow \cdots.$$
\end{proposition}

From now until the end of this section we assume that $m$ is a
multiple of $p$. From the previous proposition we clearly have
$$L(c_{p,1},h_{m,1}) \cong
V(c_{p,1},h_{m,1})/V(c_{p,1},h_{m+2p,1}).$$ Also, for every $r$
this yields a short exact sequence \be \label{ext1} 0
\longrightarrow L(c_{p,1},h_{m+(r+2)p,1}) \longrightarrow
V(c_{p,1},h_{m+rp,1})/V(c_{p,1},h_{m+(r+4)p,1}) \longrightarrow
L(c_{p,1},h_{m+rp,1}) \longrightarrow 0. \ee Similarly, dual
functor applied to (\ref{ext1}) yields \be \label{ext1dual} 0
\longrightarrow L(c_{p,1},h_{m+rp,1}) \longrightarrow
\left(V(c_{p,1},h_{m+rp,1})/V(c_{p,1},h_{m+(r+4)p,1}) \right)^*
\longrightarrow L(c_{p,1},h_{m+(r+2)p,1}) \longrightarrow 0. \ee
As in \cite{M1}, \cite{FF} or \cite{DGK} it is not hard to prove
the following result, which we state without a proof.
\begin{lemma} \label{nonzeroexts}
Let $m$ be a multiple of $p$. Extensions in (\ref{ext1}) and
(\ref{ext1dual}) are the only nontrivial non-logarithmic
extensions among two irreducible modules $L(c_{p,1},h_{m+rp,1})$
and $L(c_{p,1},h_{m+sp,1})$.
\end{lemma}

In general there could be some logarithmic extensions between
irreducible Virasoro algebra modules. For instance, the
irreducible module $L(1,1)$ admits a nontrivial self-extension,
which is logarithmic. It turns out that every non-logarithmic
self-extension of $L(c,h)$ is split exact. Indeed, suppose that
$(W,\iota,\pi)$ is a self-extension of $L(c,h)$ with the embedding
map $\iota$ and projection $\pi$. Then we consider a preimage
$w=\pi^{-1}(v_2)$ in W, where $v_2$ is the lowest weight vector in
$L(c,h)$. It is easy to see that $w$ is a lowest weight vector in
$W$ (this is not the case if $W$ is logarithmic, because $w$ can
form a Jordan block with another vector). Thus, the submodule of
$W$ generated by $w$ is a quotient of $M(c,h)$. If $v_1$ is
another lowest weight vector in $L(c,h)$, then $\iota(v_1)$
generates a copy of $L(c,h)$ in $W$. Now, $\iota(L(c,h)) \cap
U(Vir)w =0$, for if there is something nontrivial in the
intersection this would contradict to either the irreducibility of
$L(c,h)$ or the exactness. Thus, we have $\iota(L(c,h)) \oplus
U(Vir_{\leq 0})w=W$ ($W$ is generated by $w$ and $\iota(v_1)$).
Now, the kernel of $\pi$ is precisely $\iota(L(c,h)) \cong
L(c,h)$, so $U(Vir_{\leq 0})w \cong L(c,h)$, and the sequence is
split exact.

\section{Vertex operator algebra $\overline{M(1)_p}$}
\label{ver-m} In what follows the vertex operator algebra
$M(1)_{\lambda_p}$ will be denoted by $M(1)_p$, for simplicity.
This vertex operator algebra has two important vertex subalgebras.
First, $\tau$-invariant vertex subalgebra $M(1)_p^+$, where $\tau$
is the involution induced by $\tau h(-n){\bf 1}=-h(-n){\bf 1}$
(see \cite{DN} for more about $M(1)^+$ and its modules). For
$\lambda_p=0$, $M(1)^+$ is in fact a vertex operator subalgebra of
$M(1)$ (the Virasoro element is fixed by $\tau$). Another
subalgebra of interest is $\mip$, in the physics literature
usually denoted by $\mathcal{W}(2,2p-1)$ (the singlet
$\mathcal{W}$-algebra).
In this section we recall the definition of $\mip$ and some
structural results, following closely \cite{A-2003}. We will also
present a result on simplicity of $\mip $.

\subsection{Definition of $\mip$}

In what follows we shall study a subalgebra of the vertex operator
algebra $V_L$ associated to the lattice $L={\Z} \a $, $\la \a, \a
\ra =2 p$, where $p \in {\Z_{\geq 2}}$.

 Let $p \in {\N}$, $p \ge 2$. Let $L={\Z}{\a}$  be a
 lattice of rank one with nondegenerate
 ${\Z}$-
 bilinear
 form $\la \cdot, \cdot \ra$ given by
 $$ \la \a , \a \ra = 2 p .$$
 Let
${\h}={\C}\otimes_{\Z} L$. Let $\widehat{\goth{h}}$ be as in
Section 1. Extend the form $ \la \cdot, \cdot \ra $ on $L$ to
${\h}$. As in Section 1 we shall denote by $\widehat{\goth{h}}$ an
extended Heisenberg Lie algebra associated to $\goth{h}$, where
$$h=\frac{\alpha}{\sqrt{2p}},$$
and by $M(1)$ the corresponding vertex algebra, which is also a
level one $\widehat{\goth{h}}$-module. Then $M(1)$ is a vertex
subalgebra of the generalized vertex algebra $V_L$.

Consider the dual lattice $\widetilde{L}$ of $L$, so that $
\widetilde{L} = {\Z}(\frac{\alpha}{2 p})$. Let $Y$ be the vertex
operator map that defines the generalized vertex algebra structure
on  $V_{ \widetilde{L}}$. The vertex algebra $V_L$ is then a
vertex subalgebra of $V_{ \widetilde{L}}$.

As in Section 1, we shall choose the following Virasoro element in
$M(1) \subset V_L$:
$$ \omega = \frac{1 }{ 4 p } {\a} (-1) ^{2}{\bf 1} + \frac{p-1}{2p }
{\a} (-2){\bf 1}.$$ In Section 1 the corresponding vertex algebra
was denoted by $M(1)_{\lambda_p}$, but we shall write $M(1)_p$ for
simplicity.
The subalgebra of $M(1)_p$ generated by $\omega$ is isomorphic to
the simple Virasoro vertex operator algebra $L(c_{p,1},0)$ where
as before $c_{p,1}=1 - 6 \frac{(p-1) ^{2}}{p}$. Let
$$Y(\omega,z) = \sum_{n \in {\Z} } L(n) z ^{-n -2} . $$

The element $L(0)$ of the Virasoro algebra defines a $
{\N}$--gradation on $V_L$. In this article we shall consider $V_L$
as a ${\N}$--graded vertex operator algebra of rank $c_{p,1}$. As
in \cite{A-2003} define the following operators:

$$Q= e^{ \alpha} _0, \qquad \widetilde{Q} =e^{-\tfrac{1}{p} \alpha}
_0, $$ where
$$e^{\alpha}=1 \otimes e^{\alpha} \in V_{L}, \ \ \
\ e^{-\frac{\alpha}{p}}=1 \otimes e^{-\frac{\alpha}{p}} \in
V_{\tilde{L}},$$
$$Y(e^{\gamma},x)=\sum_{n \in \mathbb{Z}} e^{\gamma}_n x^{-n-1},$$
which denotes the Fourier expansion of $e^{\gamma}$. By using
results from \cite{A-2003}, we have
$$ [Q,\widetilde{Q}] =0, \ \ [L(n), Q] = [L(n), \widetilde{Q}] = 0
\ \ (n \in {\Z}).$$
Thus, the operators $Q$ and $\widetilde{Q}$ are intertwining
\footnote{This intertwining operator is not an intertwining
operator among a triple of vertex algebra modules (the two are
related though).} (or screening) operators among Virasoro algebra
modules. In fact, the Virasoro vertex operator algebra
$L(c_{p,1},0) \subset M(1)_p$ is the kernel of the screening
operator $Q$ (cf. \cite{A-2003}). Define
$$\mip = \mbox{Ker}_{M(1) } \widetilde{Q}. $$

Since $ \widetilde{Q}$ commutes with the action of the Virasoro
algebra, we have
$$L(c_{p,1},0) \subset  \mip. $$
This implies that $\mip $ is  a vertex operator subalgebra of
$M(1)_p$ in the sense of \cite{FHL} (i.e.,  $ \mip $ has the same
Virasoro element as $M(1)_p$).

The following theorem describes the structure of the vertex
operator algebra $\mip$ as a $L(c_{p,1},0)$--module.
\begin{theorem} \cite{A-2003} \label{tm-a-2003}
\item[(i)]The vertex operator algebra  $ \mip $ is a completely
reducible Virasoro algebra module and the following decomposition
holds:
\bea \mip =&& \bigoplus_{n =0} ^{\infty} U(Vir).\ u^{(n)}    \cong
\bigoplus_{n =0} ^{\infty} L(c_{p,1}, n^{2} p + n p - n),\nonumber
\eea
where \be \label{un}
 u ^{(n)} = Q^{n} e^{-n \a}.
\ee

\item[(ii)]The vertex operator algebra  $ \mip $ is generated by
$\omega$ and \be \label{H} H= Q e^{-\alpha} \ee of conformal
weight $2p-1$.

\end{theorem}
\subsection{Zhu's algebra $A(\mip)$}

We recall the definition of  Zhu's algebra for vertex operator
algebras \cite{Zh}.

Let $(V,Y, {\bf 1}, \omega)$ be a vertex operator algebra. We
shall always assume that $$V=\coprod_{ n \in {\N} } V_n, \ \
\mbox{where} \ \  V_n = \{ a \in V \ \vert \ L(0) a = n v \}. $$
For $a \in V_n$, we shall write $\wt (a) = n$.

For a homogeneous element $a \in V$ we define the bilinear maps $*
: V \otimes V \rightarrow V$, $\circ : V \otimes V \rightarrow V$
as follows: \bea
a*b &:= & \mbox{Res}_x  Y(a,x) \frac{ (1+x) ^{\wt (a) }}{x} b
= \sum_{i = 0} ^{\infty} { \wt (a) \choose i} a_{ i-1} b
, \nonumber \\
a\circ b &: = & \mbox{Res}_x  Y(a,x) \frac{ (1+x) ^{\wt (a) }}{x
^{2}} b
= \sum_{i = 0} ^{\infty} { \wt (a) \choose i} a_{ i-2} b
. \nonumber
 \eea
We extend $*$ and $\circ$ on $V \otimes V$ linearly, and denote by
$O(V)\subset V$ the linear span of elements of the form $a \circ
b$, and by $A(V)$ the quotient space $V / O(V)$. The space $A(V)$
has an associative algebra structure, the Zhu's algebra of $V$.
For instance $A(M(1)_a)$ is isomorphic to a polynomial algebra in
one variable. It is more difficult to prove
\begin{theorem}  \cite{A-2003} \label{zhu-alg}
Zhu's associative algebra $A(\mip)$ is isomorphic to the
commutative algebra \\  ${\C}[ x, y] / \la P(x,y) \ra$,
where  $\la P(x,y) \ra $ is the principal ideal generated by
\bea \label{ass-poly} P(x,y) = y ^{2} - \frac{( 4 p) ^{2 p -1} }{
(2 p -1)! ^{2}} \ ( x + \frac{(p-1) ^{2}}{4 p}) \prod_{i= 0}
^{p-2} \left( x + \frac{i}{4 p} ( 2 p - 2 -i) \right) ^{2} . \eea
\end{theorem}
The equation $P(x,y)=0$ defines a genus zero algebraic curve, with
a polynomial parametrization obtained in \cite{A-2003} (in a
special case the same formula was previously obtained in
\cite{W2}, \cite{W3}).

Let \be \label{H-F} Y(H_{-1}{\bf 1},x)=\sum_{n \in \Z} H_n
x^{-n-1}. \ee Since ${\rm wt}(H_{-1}{\bf 1})=2p-1$, it is
sometimes more convenient to shift the index in $H_n$ and work
with the generators $ H(n) = H_{ n + 2 p - 2}$.

\begin{proposition} \label{naj-tez}
Let $\mathcal{V}$ be a finite-dimensional vector space and $P \in
\mathbb{C}[x,y]$ as in (\ref{ass-poly}). Assume that
$$M=\mip \cdot \mathcal{V} =\mbox{span}_{\C} \{ a_ j v \ \vert a \in \mip, v
\in \mathcal{V} , \ j \in {\Z} \}$$ is a ${\N}$--gradable
$\mip$--module generated by the vector space $\mathcal{V}$ such
that
$$  L(m) v =\delta_{m,0} \  X  \cdot v, \ H(m) v= \delta_{m,0} \ Y \cdot  v \quad \ v  \in \mathcal{V}, m \in {\N}, $$
where $X,Y \in {\rm End}(\mathcal{V})$.
Then $P(X,Y) = 0$ as an operator on $\mathcal{V}$.
\end{proposition}
\noindent {\em Proof.} Since $M$ is $\mathbb{N}$-gradable module,
the top component of $M$ is an $A(\mip)$-module \cite{Zh}. The
rest follows from Theorem \ref{zhu-alg}. \qed

If $\mathcal{V}$ is one-dimensional, by slightly abusing language,
we say that the module $M$ from Proposition \ref{naj-tez} is a
lowest weight $\mip$--module with respect to the subalgebra ${\rm
Span}_{\mathbb{C}} \{ L(0), H(0) \}$, and the lowest weight is $(x
, y ) \in \mathbb{C}^2$. Therefore the lowest weights have to
satisfy the equation $P(x,y) = 0$.

\subsection{Simplicity of $\mip$ }

For every $n \in {\N}$ we define
$$Z_n = \mbox{Ker} _{ \mip} Q^{n+1}.$$
By using results from \cite{A-2003} we have:
$$ Z_n= \bigoplus_{i=0} ^{n} U(Vir). \  u^{(i)} \cong
\bigoplus_{i=0} ^{n} L(c_{p,1}, i ^{2} p + i p - i). $$



\begin{lemma} \label{pomoc-0}
Assume that $ n \ge 1$. There is $i \in {\N}$ such that
$$ 0 \neq H_{i} u^{(n)} \in Z_{n-1} . $$
\end{lemma}
\noindent {\em Proof.}
First we notice (see Lemma \ref{nontr-Q} below ) that \bea
\label{fel-1} Q^{2n +1} e^{-n \alpha} = 0\eea
Let $j \ge 0$. By using (\ref{un}), (\ref{H}),   (\ref{fel-1}) and
$$Q(a_j b)=(Qa)_jb+a_j (Qb) \quad (a, b \in V_L),  $$
 we obtain

\be \label{ad-fel} Q^{n}(H_j u^{(n)})=
%
%
%
\tfrac{1}{2 n +1} Q^{2 n +1} ( e^{-\a}_{j} e^{ - n \a} ) = 0. \ee
%
Therefore $H_j u ^{(n)} \in Z_{n-1}$ for every $j
\ge 0$. Assume now that $H_j u^{(n)} = 0$ for every $j \ge 0$.
Since $u^{(n)}$ is a singular vector for the Virasoro algebra, we
have that ${\C} u^{(n)}$ is the top level of the ${\N}$--graded
$\mip$--module $\mip\cdot u^{(n)}$.
Therefore $\mip\cdot u^{(n)}$ is a lowest weight $\mip$--module
with the lowest weight $(n ^2 p + n p -n, 0)$. This is a
contradiction since $P(n^{2}p + np -n,0)\ne 0$ (see Corollary
\ref{naj-tez} and Section 6 of \cite{A-2003}). So there exists
$j_1 \in {\N}$ such that $H_{j_1} u^{(n)} \ne 0$. \qed

Let us recall (cf. \cite{LL}) that a vertex operator algebra $V$
is called {\em simple} if it has no proper left ideals. In fact,
in the definition "left ideals" could be replaced by "right
ideals" or "two-sided ideals" \cite{LL}.

\begin{theorem} \label{simple}
 The vertex operator algebra $\mip$ is  simple.
\end{theorem}
\noindent {\em Proof.}
Assume that there is a (left) ideal $I$ in $\mip$ such that $0 \ne
I \ne \mip $, so there is a nonzero vector $v \in I$. Clearly, $I$
is also $L(0)$-graded. By acting with Virasoro generators on $v$
it can be easily shown that $u^{(n)} \in I$ for a certain $n \in
{\N}$. Since ${\bf 1} \notin I$, there is $n_0 \geq 1$ such that
\be \label{min-index}
 u ^{(n_0)} \in I, \quad  u^{(i)} \notin I \ \mbox{for every} \ i
< n_0. \ee
Now Lemma \ref{pomoc-0} gives that there is $j_1 \in {\N}$ such
that
$ 0 \ne H_{j_1} u^{(n_0)} \in I \cap Z _{n_0-1}$. But this
contradicts the relation (\ref{min-index}). The proof follows.
\qed

\section{Singular vectors of Feigin-Fuchs modules}

By using standard calculations in lattice vertex algebras (cf.
\cite{D}, \cite{DL},\cite{LL}, \cite{TK}) we obtain the following
important lemma:

\begin{lemma}  \label{lema-produkt-a}
In the generalized vertex algebra $V_{ \widetilde{L}}$ the
following formula holds: \bea  Y(e^{\a},x_1) \cdots
Y(e^{\a},x_{2n}) & =&
 E^{-}(-\a, x_1, \dots, x_{2n}) E^{+}(-\a, x_1, \dots, x_{2n})
\nonumber \\  & \ &\Delta_{2n} (x_1, \dots, x_{2n}) ^{2p} e ^{2n
\a}(x_1 \cdots x_{2n} ) ^{\a}
 \label{for-druga-a}
\eea
where
$$ E^{\pm}(-\a, x_1, \dots, x_{2n}) = \mbox{exp} \left( \sum_{k=1}
^{\infty}\tfrac{\alpha(\pm k)}{\pm k}( x_1 ^{\mp k} + \cdots +
x_{2n} ^{\mp k})\right)$$
and $\Delta_{2n} (x_1, \dots, x_{2n})=\prod_{i < j}(x_i-x_j)$ is
the Vandermonde determinant.
\end{lemma}

Let now $i \in \{0, \dots, p-1\}$. Define $\ga_i = \frac{i}{2p} \a
\in V_{\tilde{L}}$.

Define the operator $A= e^{\a}_{p-1-i}$. Note that in the case
$i=p-1$ we have that $A=Q$.

\begin{lemma} \label{nontr-A} We have
$$ A^{2n} e^{\ga_i -n \a}= C_n e^{\ga_i + n \a}$$ where
 the nonzero constant $C_n$ is
 $(-1) ^{n p } \frac{ (2n p)!}{p! ^{2n}}.$
\end{lemma}
\noindent {\em Proof.} First we notice that
\bea
&& A ^{2n } e^{\ga_i- n \alpha}=\mbox{Res}_{x_1} \mbox{Res}_{x_2}
\cdots \mbox{Res}_{x_{2n}} (x_1 \cdots x_{2n})^{p-i-1} \left(
Y(e^{\a},x_1) \cdots Y(e^{\a},x_{2n}) e^{\ga_i -n \a} \right)
\label{for-prva-a} .
\eea
By using Lemma \ref{lema-produkt-a}  and (\ref{for-prva-a}) we
have that
\bea
A ^{2n } e^{\ga_i- n \alpha}=C_n e^{\ga_i + n \a}
\label{for-treca-a},
\eea
where $C_n$ is the constant term in the Laurent polynomial
$$ \Delta_{2n} (x_1, \dots, x_{2n}) ^{2p} ( x_1 \cdots x_{2n} ) ^{- (2n-1)
p}. $$
Now Theorem 4.1 from \cite{Andrews} (famous Dyson's conjecture)
implies that $C_n = (-1) ^{n p } \frac{ (2n p)!}{p! ^{2n}}.$
\qed

Recall that a vector in $V_{L+\ga_i}$ is called primary if it is a
singular vector for the action of the Virasoro algebra.

Since $e ^{\ga_i-n \a} \in V_{L+\ga_i}$ is a primary vector for
every $n \in {\N}$, we have that $$Q ^{j} e ^{\ga_i-n \a}$$ is
either zero or a primary vector.

\begin{lemma} \label{nontr-Q}
Assume that $n \ge 1$.
We have: \item[(1)]
 $ Q ^{2n } e^{\ga_i- n \alpha} \ne 0$.
 \item[(2)] $Q ^{2n + j} e^{\ga_i- n \alpha} = 0$ for $j >0$.
\end{lemma}
\noindent {\em Proof.}   The assertion (1) was proven \cite{MY}
and \cite{TK}. The same assertion follows from Lemma \ref{nontr-A}
by using fact that there exists $f \in \mathcal{U}({\hh}^{+ })$,
$f \ne 0$ such that
$$ f Q^{2n} e^{\ga_i - n \a} = A^{2n} e^{\ga_i - n \a} \ne 0.$$

Note that for $j >0$, $Q ^{2n + j} e^{\ga_i- n \alpha}$ is a
singular vector of weight $h_{i+1, 2n +1}$.  But there are no
(singular) vectors in the Feigin-Fuchs module $M(1) \otimes
e^{\ga_i+ (n+j) \a}$ of this weight. \qed

\begin{remark}
\rm The non-triviality of $Q^{2n} e^{\ga_i -n \a}$ from Lemma
\ref{nontr-Q} can be also proven by using methods developed in
\cite{FF} and \cite{Fel}.
\end{remark}

\section{Irreducibility of certain $\overline{M_p(1)}$-modules}
\label{verma}

Since $\mip$ is a subalgebra of $M(1)_p$, we have that every
$M(1)_p$--module $M(1,\l)_p$ carries a natural $\mip$-module
structure. In fact, in \cite{A-2003} it was proven that every
${\N}$--graded irreducible $\mip$--module is an irreducible
subquotient of $M(1, \l)$ for certain $\l \in {\h}$.

Define  now $$\b=\frac{p-1}{2 p} \a \in \widetilde{L}.$$ In this
section we shall consider the $\mip$--module $M(1, \b) = M(1)
\otimes e^{\b}$. As we are about to see, this module is
distinguished from several point of views. We will prove that
$M(1,\b)$ is an irreducible $\mip$--module, where \be
\label{iden1} M(1,\b)=M(1,\lambda_p)_{\lambda_p} \ee is a
self-dual Virasoro module studied in Section 2.

Firstly, in parallel with Section 3, we shall view $M(1,\b)$
inside a module for the generalized lattice vertex algebra $V_L$.
For these purposes let us consider $V_L$-module $V_{L+\b}$ that
contains $M(1,\b)$. We shall now investigate the action of the
operator $Q$ on $V_{L+\b}$. Since operators $Q ^{j}$, $j \in
{\Zp}$, commute with the action of the Virasoro algebra, they are
(again) intertwining (or screening) operators among Feigin-Fuchs
modules inside the $V_L$--module $V_{L+\b}$.

Next,  we shall present a theorem describing the structure of the
 $M(1)$--module $M(1,\b)$ as a module for the Virasoro vertex
operator algebra $L(c_{p,1},0)$. The following theorem is just an
improved version of Proposition \ref{self-dual} (see \cite{FF},
\cite{Fel}).
\begin{theorem} \label{str-ff-b}
  For every $n \in {\Zp}$, the vector
$$v^{(n)}= Q^{n}
e^{\b -n \a}$$ is a non-trivial singular vector and
$$ U(Vir) v^{(n)} \cong L(c_{p,1}, h^p_n), \ \  \mbox{where} \quad
h^{p}_n= \frac{ (2 p n) ^{2} - (p-1) ^{2}}{4 p}. $$
$M(1,\b)$ is a completely reducible $L(c_{p,1},0)$--module and we
have the following decomposition
$$ M(1,\b) =
\bigoplus_{n = 0} ^{\infty} L (c_{p,1}, h^p_n).
$$
\end{theorem}

Theorem \ref{str-ff-b} immediately gives the following
result.
\begin{proposition} \label{pomoc2-b}
 We have
$$L(c_{p,1},-\tfrac{(p-1) ^{2}}{4 p}) \cong  {\rm Ker}_{M(1,\b)} Q  \ . $$
\end{proposition}

For every $n \in {\N}$ we define
$$Z^{\b}_n = \mbox{Ker} _{M(1,\b)} Q^{n+1}.$$
By using Lemma \ref{nontr-Q} and Theorem \ref{str-ff-b} we have:
$$ Z^{\b}_n= \bigoplus_{i=0} ^{n} U(Vir). \  Q^{i} e^{\b -i \a} \cong
\bigoplus_{i=0} ^{n} L(c_{p,1}, h_i^p). $$

\begin{lemma} \label{pomoc-1}
\item[(i)] There is a nonzero constant $C \in {\C}$ such that
$$ H_{j_0} v^{(n)} = C v^{(n+1)} + v', \quad v' \in Z^{\b}_n $$
where $j_0 = -2 n p + p -2$.

\item[(ii)] Assume that $ n \ge 1$.  There is $j_1 \in {\N}$ such
that
$$ H_{j_1} v^{(n)} \in Z^{\b}_{n-1}, \quad H_{j_1} v^{(n)} \ne 0 . $$
\end{lemma}
\noindent {\em Proof.}
 By using Lemma \ref{nontr-Q} we have
$$ Q^{n +1} ( H_{j_0} v^{(n)} ) = \tfrac{1}{2 n +1} Q^{2 n +2} (
e^{-\a}_{j_0} e^{\b - n \a} )= \tfrac{1}{2 n +1} Q^{2 n
+2}(e^{\beta-(n+1)\alpha}) \ne 0.
$$
Therefore $H_{j_0} v^{(n)} \in Z^{\b}_{n+1} \setminus Z^{\b}_n$.
Next we notice that
$$ L(0) H_{j_0} v^{(n)} = ( (2 p -1 ) + h_n^p - j_0 -1 ) H_{j_0} v^{(n)}
=h_{n+1}^p H_{j_0} v^{(n)}. $$
 Since $L(0) v^{(n+1)} = h_{n+1}^p v^{(n+1)}$, we conclude that
there is a constant $C$, $C \ne 0$, such that $ H_{j_0} v^{(n)} =
C v^{(n+1)} + v'$, $v' \in Z^{\b}_n $. This proves (i).

 The proof of  assertion (ii) is similar to that of Lemma
\ref{pomoc-0}.  Let $j
> j_0$. We have
$$ Q^{n } ( H_{j} v^{(n)} ) = \tfrac{1}{2 n +1} Q^{2 n +1} (
e^{-\a}_{j} e^{\b - n \a} ) = 0. $$
Therefore $H_j v ^{(n)} \in Z^{\b}_{n-1}$ for every $j > j_0$.
Assume now that $H_j v^{(n)} = 0$ for every $j \ge 0$. Since
$v^{(n)}$ is a singular vector for the Virasoro algebra, we have
that ${\C} v^{(n)}$ is the top level of the  ${\N}$--graded
$\mip$--module  $$\mip\cdot v^{(n)} = \mbox{span}_{\C} \{ a_ i
v^{(n)} \ \vert a \in \mip, \ i \in {\Z} \}.$$
Therefore $\mip\cdot v^{(n)}$ is a lowest weight $\mip$--module
with the lowest weight $(h^p_n, 0)$. This is a contradiction since
$P(h^p_n,0)\ne 0$ (see Proposition \ref{naj-tez} and Section 6 of
\cite{A-2003}). So there exists $j_1 \in {\N}$ such that $H_{j_1}
v^{(n)} \ne 0$. \qed

\begin{theorem} \label{irr-module} We have,
 $$ M(1,\b) := M(1) \otimes e ^{\b}$$ is an irreducible
 $\mip$--module ( = $\mathcal{W}(2,2p-1)$--module).
\end{theorem}
\noindent {\em Proof.} By using Lemma \ref{pomoc-1}  (i) we see
that
$v ^{(n)} \in \mip \cdot e^{\b}$ for every $n \in {\N}$. Now
Theorem \ref{str-ff-b} implies that
\bea \label{ciklicki}
&& M(1,\b) = \mip \cdot e^{\b}, \ \ \mbox{i.e.,} \ \ e^{\b} \
\mbox{is a cyclic vector.}
\eea

Assume now that there is a $\mip$--submodule $0 \ne N \subseteq
M(1,\b)$. Then  $N$ is also $L(0)$--graded. One can easily show
that $v^{(n)} \in N$ for some $n \in {\N}$. Assume that $N\ne
M(1,\b)$.  Then there is $n_0 \in {\Zp}$ such that
$$ v ^{(n_0)} \in N, \quad  v^{(i)} \notin N \ \mbox{for every} \
i < n_0.$$
Now Lemma \ref{pomoc-1}(ii) gives that there is $j_1 \in {\N}$
such that
$ 0 \ne H_{j_1} v^{(n_0)} \in N \cap Z^{\b} _{n_0-1}$. This
contradicts the minimality of $n_0$. Therefore $1 \otimes e^{\b}
\in N$, which implies that $N= M(1,\b)$. \qed


\section{Indecomposable $L(c_{p,1},0)$ and $\mip$-modules}

In the previous section we proved that $M(1,\beta)$ is an
irreducible $\mip$-module. In this part we construct a non-trivial
self-extension of $M(1,\beta)$ and describe this extension in
terms of Virasoro algebra submodules. As in \cite{M1} we start
from a two-dimensional vector space $\Omega$ and define a
$\widehat{\goth{h}}$-module
$$M(1)_p \otimes \Omega,$$
where $h(0)|_{\Omega}$ act as \be \label{jordan}
\left[\begin{array}{cc} \lambda_p & 1
\\ 0 & \lambda_p
\end{array}\right],
\ee in some basis $\{w_1,w_2 \}$ of $\Omega$. Here, we identified
${\bf 1} \otimes \Omega$ with $\Omega$. Even though the action of
$h(0)$ admits a Jordan block of size two, the module $M(1)_p
\otimes \Omega$ is still an ordinary $M(1)_p$-module (i.e., it is
diagonalizable with respect to the action of $L(0)$). Let us also
recall that the lowest conformal weight of $M(1)_p \otimes \Omega$
is \be \label{4p} \frac{-(p-1)^2}{4p}. \ee


The following result describes $M(1)_p \otimes \Omega$ in more
details.
\begin{theorem} \label{self-dual-th} Let $h(0)$ acts on $\Omega$ via (\ref{jordan}).
Then the space $M(1)_p \otimes \Omega$ is an ordinary, self-dual
(viewed as a Virasoro module), cyclic $\mip$-module. Moreover, we
have a non-split exact sequence of $\mip$-modules \be
\label{non-split} 0 \longrightarrow M(1,\beta) \longrightarrow
M(1)_p \otimes \Omega \longrightarrow M(1,\beta) \longrightarrow
0. \ee
\end{theorem}
\noindent {\em Proof:} First we observe that $M(1)_p \otimes
\mathbb{C}w_1 \cong M(1,\beta)$. Thus we have an embedding (on the
level of $L(c_{p,1},0)$ and $\mip$-modules) $M(1,\beta)
\hookrightarrow M(1)_p \otimes \Omega$. Also, from
$h(0)w_2=\lambda_p w_2 +w_1$ it is clear that $M(1)_p \otimes
\Omega/M(1,\beta) \cong M(1,\beta)$. Thus we have an exact
sequence (\ref{non-split}).

To show that $M(1)_p \otimes \Omega$ is a self-dual Virasoro
algebra module, it suffices to show that the operator $L(0)$
acting on $\Omega^*$ is similar to the operator $L(0)$ acting now
on $\Omega$. This is a consequence of the following more general
fact that applies for any pair of modules: Let ${\rm
dim}(\Omega_1)={\rm dim}(\Omega_2)$. Then $M(1)_p \otimes
\Omega_1$ and $M(1)_p \otimes \Omega_2$ are isomorphic as
$M(1)_p$-modules if and only if $h(0)|_{\Omega_1}$ is similar to
$h(0)|_{\Omega_2}$. Now, let us focus on the contragradient module
$(M(1)_p \otimes \Omega)^*$. By using the same argument as in
Lemma \ref{dual-iso} it follows that
$$(M(1)_{p} \otimes \Omega)^* \cong M(1)_{p}
\otimes \Omega^*,$$ where $h(0)|_{\Omega^*}$ is represented by
$$\left[\begin{array}{cc} \lambda_p & 0 \\ -1
& \lambda_p \end{array}\right].$$ But this module is isomorphic to
$M(1)_{p} \otimes \Omega$, so our module is self-dual. Now, we
prove that the extension in (\ref{non-split}) is non-split. For
these purposes we compute the action of $H(0)$ on $\Omega$. From
\cite{A1}, \cite{A-2003} we have
$$H(0) \cdot w_1=\frac{\alpha(\alpha-1) \cdots (\alpha-2p+2)}{(2p-1)!} \cdot
w_1,$$ where $\alpha$ acts as
$$\left[\begin{array}{cc} p-1 & \sqrt{2p} \\ 0 & p-1 \end{array} \right].$$
Now if we combine the previous two formulas we get that $H(0)$
acts (in the same basis) via the nilpotent Jordan block
$$\left[ \begin{array}{cc} 0 & \nu_p \\ 0 & 0 \end{array} \right].$$
Finally, since $H(0) w_2= \nu_p w_1$, for some $\nu_p \neq 0$ and
$M(1,\beta)$ is irreducible (and hence cyclic) it is clear that
$M(1)_p \otimes \Omega$ is also cyclic. \qed

\begin{remark}
\rm The $\mip$-module $M(1)_p \otimes \Omega$ can be also
constructed by using the Zhu's theory (cf. \cite{Zh}) starting
from a two-dimensional $A(\mip)$-module $\Omega$.
\end{remark}

Here is an interesting consequence of formula (\ref{ass-poly}),
which also indicates why $M(1,\beta)$ is indeed a special
$\mip$-module.

\begin{proposition}
There are no logarithmic self-extensions of $M(1,\beta)$.
\end{proposition}
\noindent {\em Proof.} Suppose that
$$0 \longrightarrow M(1,\beta) \longrightarrow M \longrightarrow
M(1,\beta) \longrightarrow 0,$$ for some logarithmic module $M$.
Since $M(1,\beta)$ is irreducible, $M$ is generated by $2$ vectors
$w_1$ and $w_2$ of generalized conformal weight $-(p-1)^2/4p$.
Thus, we may assume that $w_1,w_2$ form a Jordan block with
respect to $L(0)$, that is $L(0)w_1=\frac{-(p-1)^2}{4p}w_1$ and
$L(0)w_2=\frac{-(p-1)^2}{4p}w_2+w_1$. Now, for every
$\mathbb{N}$-gradable $\mip$-module $M=\oplus_{n \in \mathbb{N}}
M_n$, the top component $M_0$ has a natural $A(\mip)$-module
structure. By Proposition \ref{naj-tez}, $P(L(0),H(0))=0$ as an
operator on acting on $M_0$. The equation $P(L(0),H(0))=0$ implies
that $H^2(0) \cdot w_1=0$ and $H^2(0) \cdot w_2=a w_1$, for some
$a \neq 0$, which depends on $p$. But there is no linear operator
with these properties. \qed

\section{Virasoro algebra structure of $M(1)_p \otimes \Omega$}

In this part we describe $M(1)_p \otimes \Omega$ as a Virasoro
algebra module. The embedding structure of $M(1)_p \otimes \Omega$
is similar to the case $c=1$ studied in \cite{M1}. As in the
previous section we shall denote by $\{w_1,w_2 \}$ a basis of
$\Omega$ with the action of $L(0)$ and $H(0)$ computed as above
via (\ref{jordan}). In what follows we shall use the following
graphical notation: a singular vector in $M(1)_p \otimes \Omega$
will be denoted by $\bullet$, and a cosingular vector (i.e., a
vector that becomes singular in the quotient generated by all
singular vectors) will be denoted by $\diamond$. Then we have the
following theorem

\begin{theorem} As a Virasoro algebra module,
$M(1)_p \otimes \Omega$ is generated by a sequence of singular
vectors $v^{(n)}, n \geq 0$ and a sequence of cosingular vectors
$v^{2,n}, n \geq 0$, of conformal weight $h_n^p$, $n \geq 0$, as
on the following diagram: \be 
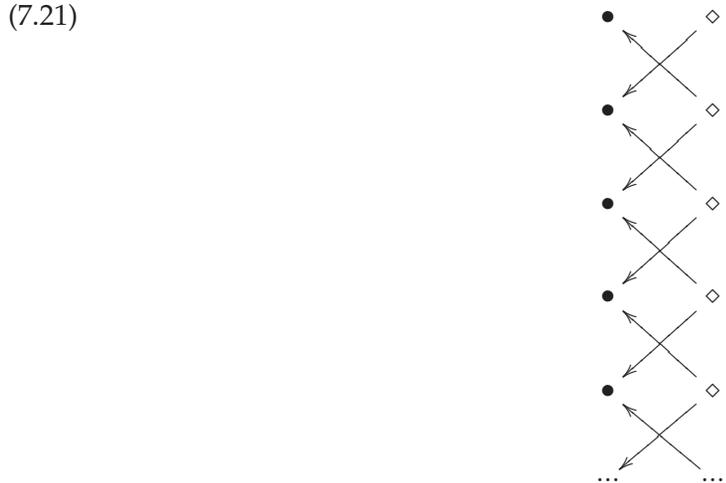
\label{embed} \xymatrix{ & {\bullet}
&
{} \ar[ld] \diamond \\
{} & {\bullet}  & {\diamond} \ar[lu] \ar[ld]
\\ {} & {\bullet} & {\diamond} \ar[lu] \ar[ld]  \\ {} &
{\bullet}
& {} \ar[lu] \ar[ld] \diamond  \\
{} & {\bullet} & {\diamond} \ar[lu] \ar[ld]
\\ {} & .. . &  ...\ar[lu] }
\ee where singular vectors are denoted by $\bullet$ and cosingular
vectors with $\diamond$. The arrows have usual meaning (i.e., an
arrow pointing from $\diamond$ to $\bullet$ indicates that
$\bullet$ is in the submodule generated by $\diamond$).
\end{theorem}
\noindent {\em Proof.} The proof is similar as in the $c=1$ case
examined in \cite{M1} so we omit some details. Let $w_1$ and $w_2$
be as in Section 6, so that
$$h(0) \cdot w_1=\lambda_p w_1, \ \ h(0) \cdot w_2=\lambda_p
w_2+w_1.$$ By using Theorem \ref{self-dual-th} it is clear that
$M(1,\beta) \hookrightarrow M(1)_p \otimes \Omega$ and a sequence
of lowest weight vectors $v^{(n)}$, described in the left column
(\ref{embed}), corresponds to the sequence of lowest weight
vectors in $M(1,\beta)$. For the cosingular vectors we choose any
sequence $v^{2,n}$ satisfying $h(0) \cdot v^{2,n}=\lambda_p
v^{2,n}+v^{(n)}$, so that $v^{2,n}$ become a (nonzero) singular
vectors in the quotient $M(1)_p \otimes \Omega/M(1,\beta)$. These
two sequences generate the whole module $M(1)_p \otimes \Omega$.
As in \cite{M1}, by an application of Lemma \ref{nonzeroexts}, we
see that the set of arrows exiting from $v^{2,n}$ is a subset of
arrows already displayed on the diagram (\ref{embed}). So it
remains to show that from every $v^{2,n}$, $n \geq 1$ there are
precisely two arrows, one pointing to $v^{(n-1)}$ and another
pointing to $v^{(n+1)}$ (from $v^{2,0}$ there should be only one
arrow pointing to $v^{(1)}$). Firstly, we show that from every
subsingular vector $v^{2,n}$, $n \geq 1$ there is an arrow
pointing up to the singular vector $v^{(n-1)}$. In order to see
this consider
$$L(m) v^{2,n}, \ \ m \geq 1.$$
We claim that for every $n \geq 1$ there exists $m \geq 1$ such
that $L(m) v^{2,n} \neq 0$. Suppose that $L(m)v^{2,n}=0$, for all
$m >0$.  This amounts to $h(m) v^{(n)}=0$, $m>0$, which is
impossible for $n \geq 1$. Since
$${\rm wt}(L(m)v^{2,n}) < {\rm wt}(v^{2,n}), \ m \geq 1,$$
then $L(m)v^{2,n} \in U(Vir)v^{(n-1)}$. It remains to show that
from every subsingular vector $v^{2,n}$ there is an arrow pointing
down to the lowest weight vector $v^{(n+1)}$. For these purposes
let us consider $(M(1)_p \otimes \Omega)^*$, the dual of $M(1)_p
\otimes \Omega$. The duality functor reverses the roles of
singular and cosingular vectors and reverses the orientation of
arrows. More precisely, there is an isomorphism from $M(1)_p
\otimes \Omega$ to its dual that maps singular vector $v^{(n)}$ to
a cosingular vector $w^{2,n}$ and the cosingular vector $v^{2,n}$
to a singular vector $w^{(n)}$, such that $w^{(n)}$ and $w^{2,n}$
form a Jordan block with respect to $h(0)$ (i.e., $h(0) \cdot
w^{(n)}=\lambda_p w^{(n)}$, $h(0) \cdot w^{2,n}=\lambda_p
w^{2,n}+w^{(n)}$).
Thus, there will be a sequence of arrows pointing {\em down} from
each cosingular vector $w^{2,n}$ to the singular vector
$w^{(n+1)}$. By using the same argument as before, from each
cosingular vector $w^{2,n}$ in $(M(1)_p \otimes \Omega)^*$ there
will be an arrow pointing up to $w^{(n-1)}$. Finally, $M(1)_p
\otimes \Omega$ is self-dual, so the proof follows. \qed

\section{A realization of self-dual $\mip$--modules }

From Theorem \ref{zhu-alg}  follows that the irreducible self-dual
$\mip$--modules have lowest weight $(x,0)$ where $ x=
-\frac{i(2p-2-i)}{4p }$, $i= 0, \dots, p-1$. (Extremal) self-dual
modules with lowest weights $(0,0)$ and
$(-\frac{(p-1)^{2}}{4p},0)$ were constructed in Sections
\ref{ver-m} and \ref{verma}.

As before for  $i \in \{0, \dots, p-2\}$, we  define $\ga_i = \frac{i}{2p}\a$, and consider $V_L$--module $V_{L+\ga_i}$.

Recall that $h_{m,n}:=\frac{(m-np)^2-(p-1)^2}{4p}.$

By using Lemma \ref{nontr-Q} and  the structure theory of
Feigin-Fuchs modules \cite{FF}  we get the following theorem.

\begin{theorem} \label{str-ff}
Assume that $i \in \{0, \dots, p-2\}$.
 \item[(i)] The Feigin-Fuchs module $M(1, \ga_i)$, is generated by the family of
singular and cosingular vectors $ \widetilde{Sing}_{i}  \bigcup
\widetilde{CSing}_{i}$, where
$$  \widetilde{Sing}_{i} =  \{ u_i ^{(n)} \ \vert \ n \in {\N}\}; \
\
  \widetilde{CSing}_{i} =  \{ w_i ^ {(n)} \ \vert \ n \in {\Zp} \}. $$
These vectors satisfy the following relations:
\bea
&&  u_i ^{(n)}= Q ^{n} e ^{\ga_i -n \alpha} , \ \ Q ^{n} w_i ^{(n)} = e ^{\ga_i + n \a},
\nonumber \\
&& U(Vir) u_i ^{(n)} \cong L(c_{p,1}, h_{i+1, 2n +1}). \nonumber
\eea
\item[(ii)] The submodule generated by vectors $u_i ^{(n)}, n \in {\N}$ is
isomorphic to $$ \overline{M(1,\ga_i}) \cong \bigoplus_{n =0 } ^{\infty}
L(c_{p,1}, h_{i+1, 2n+1}).$$

\item[(iii)]The quotient module  is isomorphic to
$$M(1,\ga_i)/  \overline{M(1,\ga_i}) \cong \bigoplus_{n =1 } ^{\infty}
L(c_{p,1}, h_{i+1, -2n +1}). $$
\end{theorem}

\begin{theorem} \label{irreducible-ga-i}
   Assume that $i \in \{0, \dots, p-2\}$. Then $\overline{M(1,\ga_i)}$ is an irreducible $\mip$--module. Moreover, $\overline{M(1,\ga_i)}$
 is an irreducible, self-dual $\mip$--module with the lowest weight  $(-\frac{i(2p-2-i)}{4p }, 0)$.
\end{theorem}
{\em Proof.} By using Frenkel-Zhu's formula \cite{FZ} and the
methods developed in \cite{M0} one can prove that the space of
intertwining operators
$$I \ { L(c_{p,1},h) \choose  L(c_{p,1},2 p -1) \ \  L(c_{p,1},
h_{i+1,2n+1})}$$ is non-trivial if and only if $h=h_{i+1,2n-1}$,
$h=h_{i+1,2n+1}$ or $h=h_{i+1,2n+3}$, for $n \geq 1$. Since the
multiplicities of these fusion rules are always one, we may write
formally:
\bea  &&
L(c_{p,1},2 p -1) \times L(c_{p,1}, h_{i+1,2n+1} ) =\nonumber \\
 && L(c_{p,1}, h_{i+1,2n-1} ) \oplus L(c_{p,1}, h_{i+1,2n+ 1} ) \oplus L(c_{p,1}, h_{i+1,2n+3} ) \quad (n \ge 1).\label{f-r}\eea
The same results has been known by physicists (cf. \cite{F1}). The
fusion rules (\ref{f-r}) implies that
$$H_j \overline{M(1,\ga_i)} \subset \overline{M(1,\ga_i)} , \quad
\mbox{for every} j \in {\Z}.$$ Since $\mip$ is generated by
$\omega$ and $H$ we have that $\overline{M(1,\ga_i)}$ is an
$\mip$--module. By using a completely analogous proof to those of
Theorem \ref{simple} and Theorem \ref{irr-module} we get the
$\overline{M(1,\ga_i)}$ is an irreducible $\mip$--module. \qed

\section{Some logarithmic $\mip$-modules}

In this part we study certain logarithmic $M(1)_p$-modules. First,
we consider $M(1)_p$-module $M(1)_p \otimes \Omega_0$, where on
$\Omega_0$ we have $h(0)^2=0$ and $h(0) \neq 0$. Then we have
$$0 \lar M(1)_p \lar M(1)_p \otimes \Omega_0 \lar M(1)_p \lar 0.$$
Furthermore, $L(0)^2=0$ and $L(0) \neq 0$ on $\Omega_0$, so
$L(0)|_{\Omega_0}$ is a nilpotent Jordan block of size two. The
module $M(1)_p \otimes \Omega_0$ is also an $\mip$-module, where
on the two-dimensional top level $\Omega$, the generator $H(0)$
acts via a nonzero nilpotent Jordan block (this matrix can be
easily computed by using Proposition 5.2 in \cite{A-2003} or
relation (\ref{ass-poly})). In order to describe the structure of
$M(1)_p \otimes \Omega_0$ let us first examine $M(1)_p$, viewed as
an $\mip$-module. As before, we shall view $M(1)_p$ embedded
inside the lattice VOA $V_L$, which further sits inside the
generalized vertex algebra $V_{\tilde{L}}$. Now, the Virasoro
algebra dual of $M(1)_p$ is also contained in $V_{\tilde{L}}$ and
is isomorphic to $M(1)_p \otimes e^{2 \beta}=M(1,2 \beta)$, by
Lemma \ref{dual-iso}. Consider a two-dimensional space $\Omega_1$
such that $(h(0)-2 \lambda_p)^2=0$ and $h(0) \neq 2 \lambda_p $,
as operators on $\Omega_1$. Then we have an exact sequence
$$0 \lar M(1,2 \beta) \lar M(1)_p \otimes \Omega_1 \lar M(1,2
\beta) \lar 0,$$ where $M(1)_p \otimes \Omega_1$ is a logarithmic
module $\mip$-module (keep in mind that
$h(0)=\frac{\alpha}{\sqrt{2p}}$ and
$\beta=\frac{p-1}{\sqrt{2p}}$).

We observe an exact sequence of $\mip$-modules
$$0 \lar \mip \lar M(1)_p \lar W \lar 0,$$
where $W \cong M(1)_p/\mip$ is an $\mip$-module. Recall that the
screening $\tilde{Q}=e^{-\alpha/p}_0$ acts as a derivation
operator on $V_L$, meaning that
$$\tilde{Q}(u_j v)=(\tilde{Q}u)_j v+u_j \tilde{Q} v, \ \ u,v \in V_L.$$
Also, by definition (cf. \cite{A-2003})
$${\rm Ker}|_{M(1)_p} \tilde{Q}=\mip.$$
In view of that, the map
$$\tilde{Q} : M(1)_p \longrightarrow M(1)_p \otimes
e^{-\alpha/p}$$ is actually a homomorphism of $\mip$-modules, so
we have an isomorphism
$${\rm Im} \ \tilde{Q} \cong W.$$
So $W$ is actually an $\mip$-submodule of $M(1)_p \otimes
e^{-\alpha/p}$. This situation is depicted via the following
diagram
$$
\xymatrix{ & & \bullet  &  & \\ & \bullet \ar@{.}[ur] &
\ar@/^/[l]^{\tilde{Q}}
\diamond \ar[u] \ar[d] \ar@/_/[r]_Q & \ar@{.}[ul] \diamond \ar[d] & \\
\diamond \ar[d] \ar@{.}[ur] & \diamond \ar[u] \ar[d] & \bullet &
\bullet  & \ar@{.}[ul] \bullet \\ ... & ... & ... \ar[u] & ...
\ar[u] & ... \ar[u] }
$$
The middle column represents the vacuum module $M(1)_p$. The
column to the right (resp. left) of $M(1)_p$ represents $M(1)_p
\otimes e^{\alpha}$ (resp. $M(1)_p \otimes e^{-\alpha/p}$).
Observe also that
$$W \cong \bigoplus_{n=1}^\infty L(c_{p,1},n^2p-np+n).$$
In addition $W$ is a cyclic $\mip$-module of lowest weight
$(L(0),H(0))=(1,-2p)$.

Now, $M(1,2 \beta)$ enjoys similar properties. Firstly, we observe
a sequence
$$0 \lar W^* \lar M(1,2 \beta) \lar \mip \lar 0,$$
where we used self-duality of $\mip$. As a Virasoro algebra module
$$W^* \cong \bigoplus_{n=1}^\infty L(c_{p,1},n^2p-np+n).$$
However, $W$ is not isomorphic to $W^*$ (the lowest weight of
$W^*$ is $(1,2p)$). On the other hand, both $\mip$ and
$M(1,\beta)$ are self-dual as $\mip$-modules.

\begin{remark} \rm There are some obvious questions that could be
pursued now regarding $W$ and other related modules (e.g.,
irreducibility). We shall address these issues in our future
publications.
\end{remark}

\section{Graded dimensions of self-dual $\mip$-modules}

In this section we discuss graded dimensions of self-dual
$\mip$-modules appearing in the paper. As usual the graded
dimensions are defined by using the formula
$${\rm ch}_{M}(\tau)={\rm tr}|_{M} q^{L(0)-c_{p,1}/24}, \ q=e^{2 \pi i \tau}.$$
Firstly, the characters of $M(1,\beta)$ and $M(1)_p \otimes
\Omega$ are easily computed: \be {\rm
ch}_{M(1,\beta)}(\tau)=\frac{1}{\eta(\tau)}, \ee \be {\rm
ch}_{M(1)_p \otimes \Omega}(\tau)=\frac{2}{\eta(\tau)}. \ee It is
also not hard to see, by using Proposition \ref{VermaVir} and
Theorem \ref{irreducible-ga-i}, that  \be {\rm
ch}_{\overline{M(1,\gamma_i)}}(\tau)=q^{(i-p+1)^2/4p} \cdot
\frac{\displaystyle{\sum_{n=0}^\infty q^{n(p(n+1)-i-1)}-
\sum_{n=1}^\infty q^{n(pn-p+i+1)}}}{\eta(\tau)}, \ee where
$i=0,...,p-2$. The previous formula can be rewritten in a more
compact way as \be \label{charself} {\rm
ch}_{\overline{M(1,\gamma_i)}}(\tau)=\frac{\displaystyle{\sum_{n
\in \mathbb{Z}} {\rm sgn}(n)
q^{\frac{((2n+1)p-i-1)^2}{4p}}}}{\eta(\tau)}, \  \ i=0,..,p-2, \ee
where ${\rm sgn}$ is the sign function (where ${\rm sgn}(0)=1$.)
Just for the record notice also that
$${\rm
ch}_{M(1)_p \otimes \Omega_0}(\tau)={\rm ch}_{M(1)_p \otimes
\Omega_1}(\tau)=\frac{q^{(p-1)^2/4p}(1+2 \pi i
\tau)}{\eta(\tau)}.$$

The numerators appearing in (\ref{charself}) are certain
$\theta$-like constants. Their modular properties were studied in
\cite{F1}.


\section{Construction of hidden logarithmic intertwining operators for $\mathcal{W}(2,2p-1)$-algebras}

In this part we apply results from \cite{M1} and the previous
section in the case of $\lambda=\lambda_p$. Results from \cite{M1}
provide us with a hidden logarithmic intertwining operator
$$\mathcal{Y} \in I \ { M(1)_p \otimes \Omega_1 \choose M(1)_p
\otimes \Omega \ \ M(1,\beta)},$$ where $M(1)_p \otimes \Omega_1$
is the logarithmic $\mip$-module described in the previous
section. The explicit formulas for $\mathcal{Y}$ are given in
Theorem 7.4, \cite{M1}. Here we give explicit formulas for
$\mathcal{Y}(u,x)$ in two special cases: $u=w_i \in M(1)_p \otimes
\Omega$, $i=1,2$.

\begin{corollary} \label{hidden}
Let $w_1$, $w_2$ and $\mathcal{Y}$ as above. Then we have
 \bea
\mathcal{Y}(w_1,x)&=&E^-(\lambda_p,x)E^+(\lambda_p,x)T_{\Omega,
\beta}^{\Omega_1}(w_1)(1+ \lambda_p {\rm log}(x) h_n(0))
x^{\lambda_p h_s(0)} \nn
\mathcal{Y}(w_2,x)&=& \int^- h(x)
E^-(\lambda_p,x)E^+(\lambda_p,x)T_{\Omega,
\beta}^{\Omega_1}(w_1)(1+ \lambda_p {\rm log}(x) h_n(0))
x^{\lambda_p h_s(0)} \nn
&& + E^-(\lambda_p,x)E^+(\lambda_p,x)T_{\Omega,
\beta}^{\Omega_1}(w_1)(1+ \lambda_p {\rm log}(x) h_n(0))
x^{\lambda_p h_s(0)} \int^+ h(x) \nn
&&+E^-(\lambda_p,x)E^+(\lambda_p,x)T_{\Omega,
\beta}^{\Omega_1}(w_2)(1+ \lambda_p {\rm log}(x) h_n(0))
x^{\lambda_p h_s(0)}, \eea
where $h_n(0)$ is the nilpotent part of $h(0)$, $T_{\Omega,
\beta}^{\Omega_1} \in Hom(\Omega, Hom(\mathbb{C}_{\beta},
\Omega_1))$ is the operator that corresponds to the obvious
$\goth{h}$-isomorphism between $\Omega \otimes \mathbb{C}_{\beta}$
and $\Omega_1$, and \bea \int^{+} h(x)&=&h(0) \log(x) + \sum_{m
>0} \frac{h(m) x^{-m}}{-m}, \nn \int^{-} h(x)&=&\sum_{m < 0}
\frac{h(m) x^{-m}}{-m}. \nonumber \eea

\end{corollary}

\section{Conclusion and future work}

As we indicated in the introduction our future goal(s) are in the
direction of understanding the triplet vertex operator algebra
$\mathcal{W}(2,(2p-1)^3)$ \cite{F3}, \cite{Ga},
\cite{GK1},\cite{GK2} and its irreducible modules. This will
require some modifications of the present paper since the triplet
algebra $\mathcal{W}(2,(2p-1)^3)$ is not contained inside
$M(1)_p$, but rather inside $V_{L}$. Another interesting
direction, which in our opinion does not have a satisfactory
explanation, is to probe our method for certain logarithmic
$A_1^{(1)}$-modules at admissible level studied by Gaberdiel
\cite{Ga}.

\end{document}